\documentclass[12pt,a4paper]{article}
\usepackage{a4,amsmath,amssymb,amsfonts,amsthm}
\usepackage{amsrefs}
\usepackage[all]{xy}
\usepackage[small,nohug,heads=vee]{diagrams}

\begin{document}
\bibliographystyle{plain}
 

\def\mR{\M{R}}           
\def\mZ{\M{Z}}           
\def\mN{\M{N}}           
\def\mQ{\M{Q}}       
\def\mC{\M{C}}  
\def\mG{\M{G}}



\def\Spec{{\rm Spec}}
\def\rg{{\rm rg}}
\def\Hom{{\rm Hom}}
\def\Aut{{\rm Aut}}
 \def\Tr{{\rm Tr}}
 \def\Exp{{\rm Exp}}
 \def\Gal{{\rm Gal}}
 \def\End{{\rm End}}
 \def\det{{{\rm det}}}
 \def\Td{{\rm Td}}
 \def\ch{{\rm ch}}
 \def\che{{\rm ch}_{\rm eq}}
  \def\Spec{{\rm Spec}}
\def\Id{{\rm Id}}
\def\Zar{{\rm Zar}}
\def\Supp{{\rm Supp}}
\def\eq{{\rm eq}}
\def\Ann{{\rm Ann}}
\def\LT{{\rm LT}}
\def\Pic{{\rm Pic}}
\def\rg{{\rm rg}}
\def\et{{\rm et}}
\def\sep{{\rm s}}
\def\ppcm{{\rm ppcm}}
\def\ord{{\rm ord}}
\def\Gr{{\rm Gr}}
\def\ker{{\rm ker}}
\def\rk{{\rm rk}}
\def\Stab{{\rm Stab}}
\def\im{{\rm im}}
\def\Sm{{\rm Sm}}
\def\red{{\rm red}}
\def\Frob{{\rm Frob}}
\def\Ver{{\rm Ver}}


\def\beginProof{\par{\bf Proof }}
 \def\endProof{${\qed}$\par\smallskip}
 \def\pr{^{\prime}}
 \def\prpr{^{\prime\prime}}
 \def\mtr#1{\overline{#1}}
 \def\ra{\rightarrow}
 \def\mfp{{\mathfrak p}}
 \def\mfm{{\mathfrak m}}
 
 \def\mQ{{\Bbb Q}}
 \def\mR{{\Bbb R}}
 \def\mZ{{\Bbb Z}}
 \def\mC{{\Bbb C}}
 \def\mN{{\Bbb N}}
 \def\mF{{\Bbb F}}
 \def\mA{{\Bbb A}}
  \def\mG{{\Bbb G}}
 \def\CI{{\cal I}}
 \def\CA{{\cal A}}
 \def\CE{{\cal E}}
 \def\CJ{{\cal J}}
 \def\CH{{\cal H}}
 \def\CO{{\cal O}}
 \def\CA{{\cal A}}
 \def\CB{{\cal B}}
 \def\CC{{\cal C}}
 \def\CK{{\cal K}}
 \def\CL{{\cal L}}
 \def\CI{{\cal I}}
 \def\CM{{\cal M}}
\def\CP{{\cal P}}
 \def\CZ{{\cal Z}}
\def\CR{{\cal R}}
\def\CG{{\cal G}}
\def\CX{{\cal X}}
\def\CY{{\cal Y}}
\def\CV{{\cal V}}
\def\CW{{\cal W}}
 \def\wt#1{{\widetilde{#1}}}
 \def\mod{{\rm mod\ }}
 \def\refeq#1{(\ref{#1})}
 \def\blb{{\big(}}
 \def\brb{{\big)}}
\def\mc{{{\mathfrak c}}}
\def\mcpr{{{\mathfrak c}'}}
\def\mcprpr{{{\mathfrak c}''}}
\def\ss{{\rm ss}}
\def\parf{{\rm parf}}
\def\P1{{{\bf P}^1}}
\def\cod{{\rm cod}}
\def\pr{\prime}
\def\prpr{\prime\prime}
\def\ss{\scriptstyle}
\def\OX{{ {\cal O}_X}}
\def\mpartial{{\mtr{\partial}}}
\def\inv{{\rm inv}}
\def\indlim{\underrightarrow{\lim}}
\def\prolim{\underleftarrow{\lim}}
\def\pprolim{'\prolim'}
\def\Pro{{\rm Pro}}
\def\Ind{{\rm Ind}}
\def\Ens{{\rm Ens}}
\def\without{\backslash}
\def\pbdb{{\Pro_b\ D^-_c}}
\def\qc{{\rm qc}}
\def\Com{{\rm Com}}
\def\an{{\rm an}}
\def\gfield{{\rm\bf k}}
\def\s{{\rm s}}
\def\dR{{\rm dR}}
\def\ari#1{\widehat{#1}}
\def\ul#1{\underline{#1}}
\def\sul#1{\underline{\scriptsize #1}}
\def\mou{{\mathfrak u}}
\def\ich{\mathfrak{ch}}
\def\cl{{\rm cl}}
\def\K{{\rm K}}
\def\R{{\rm R}}
\def\F{{\rm F}}
\def\L{{\rm L}}
\def\pgcd{{\rm pgcd}}
\def\rc{{\rm c}}
\def\N{{\rm N}}
\def\E{{\rm E}}
\def\H{{\rm H}}
\def\CHOW{{\rm CH}}
\def\A{{\rm A}}
\def\d{{\rm d}}
\def\Res{{\rm  Res}}
\def\GL{{\rm GL}}
\def\Alb{{\rm Alb}}
\def\alb{{\rm alb}}
\def\Hdg{{\rm Hdg}}
\def\Num{{\rm Num}}
\def\Irr{{\rm Irr}}
\def\Frac{{\rm Frac}}
\def\Sym{{\rm Sym}}
\def\TV{\rm TV}
\def\indlim{\underrightarrow{\lim}}
\def\prolim{\underleftarrow{\lim}}
\def\Ver{{\rm Ver}}
\def\hn{{\rm hn}}
\def\min{{\rm min}}
\def\max{{\rm max}}
\def\Div{{\rm Div}}
\def\sm{{\rm sm}}


\def\RHom{{\rm RHom}}
\def\rRHom{{\mathcal RHom}}
\def\rHom{{\mathcal Hom}}
\def\dotimes{{\overline{\otimes}}} 
\def\Ext{{\rm Ext}}
\def\rExt{{\mathcal Ext}}
\def\Tor{{\rm Tor}}
\def\rTor{{\mathcal Tor}}
\def\SP{{\mathfrak S}}
\def\perf{{\rm perf}}

\def\H{{\rm H}}
\def\D{{\rm D}}
\def\Del{{\mathfrak D}}

\def\sh{{\rm sh}}

 \newtheorem{theor}{Theorem}[section]
 \newtheorem{prop}[theor]{Proposition}
 \newtheorem{propdef}[theor]{Proposition-Definition}
 \newtheorem{sublemma}[theor]{sublemma}
 \newtheorem{cor}[theor]{Corollary}
 \newtheorem{lemma}[theor]{Lemma}
 \newtheorem{sublem}[theor]{sub-lemma}
 \newtheorem{defin}[theor]{Definition}
 \newtheorem{conj}[theor]{Conjecture}
 \newtheorem{quest}[theor]{Question}

 \parindent=0pt
 \parskip=5pt

 \author{
 Damian R\"OSSLER\footnote{Institut de Math\'ematiques, 
Equipe Emile Picard, 
Universit\'e Paul Sabatier, 
118 Route de Narbonne,  
31062 Toulouse cedex 9, 
FRANCE, E-mail: rossler@math.univ-toulouse.fr}}
 \title{On the group of purely inseparable points of an abelian variety defined over a function field of positive characteristic}
\maketitle

\begin{abstract}
Let $K$ be the function field of a smooth and proper curve $S$ over an algebraically closed field $k$ of characteristic $p>0$. Let 
$A$ be an ordinary abelian variety over $K$. Suppose that the N\'eron model $\CA$ of $A$ over $S$ has a closed fibre $\CA_s$, which is 
an abelian variety of $p$-rank $0$. We show that under these assumptions the group 
$A(K^\perf)/\Tr_{K|k}(A)(k)$ is finitely generated. Here $K^\perf=K^{p^{-\infty}}$ is the maximal purely inseparable extension of $K$. This result implies that in some circumstances, the "full" Mordell-Lang conjecture, as well as a conjecture of Esnault and Langer, are verified. 
\end{abstract}

\section{Introduction}

Let $k$ be an algebraically closed  field of characteristic $p>0$ and let 
$S$ be a connected, smooth and proper curve over $k$. Let $K:=\kappa(S)$ be its function field. 

If $V/S$ is a locally free coherent sheaf on $S$, we denote by 
$$
0=V_0\subseteq V_1\subseteq V_2\subseteq\dots\subseteq V_{\hn(V)}=V
$$
the Harder-Narasimhan filtration of $V$. We write as usual
$$\deg(\ast):=\deg({\rm c}_1(\ast)),\ \mu(\ast):=\deg(\ast)/\rk(\ast)$$
and 
$$\mu_\min(V):=\mu(V/V_{\hn(V)-1}),\ \mu_\max(V):=\mu(V_{1}).$$ 
See  \cite[chap. 5]{Brenner-Herzog-Villamayor-Three} (for instance) for the definition of the Harder-Narashimha filtration and for the notion of semistable sheaf, which underlies it.

A locally free sheaf 
$V$ on $S$ is said to be strongly semistable if $F^{r,*}_S(V)$ is semistable for all $r\in\mN$. Here $F_S$ is the absolute Frobenius endomorphism of $S.$  A. Langer proved in \cite[Th. 2.7, p. 259]{Langer-Semistable} that  
there is an $n_0=n_0(V)\in\mN$ such that the quotients of the Harder-Narasimhan filtration of $F^{n_0,*}_S(V)$ are 
all strongly semistable. This shows in particular that the following definitions : 
$$
\bar\mu_\min(V):=\lim_{l\to\infty}\mu_\min(F^{l*}_S(V))/p^l
$$
and 
$$
\bar\mu_\max(V):=\lim_{l\to\infty}\mu_\max(F^{l*}_S(V))/p^l
$$
make sense. 

With these definitions in hand, we are now in a position to formulate the results that we are going 
to prove in the present text. 

Let $\pi:\CA\to S$ be a smooth commutative group scheme and let $A:=\CA_K$ be the generic fibre of $\CA$. 
Let $\epsilon:S\to\CA$ be the zero-section and let $\omega:=\epsilon^*(\Omega^1_{\CA/S})$ be the Hodge 
bundle of $\CA$ over $S$. 

Fix an algebraic closure $\bar K$ of $K$. For any $\ell\in\mN$, let 
$$K^{p^{-\ell}}:=\{x\in\bar K| x^{p^l}\in K\}, $$ which is a field. We may then define the field 
$$K^\perf=K^{p^{-\infty}}=\cup_{\ell\in\mN}K^{p^{-\ell}},$$
which is often called the {\it perfection} of $K$. 

\begin{theor}
Suppose that  $\CA/S$ is semiabelian and that $A$ is a principally polarized abelian variety.  Suppose that the vector bundle $\omega$ is ample. Then there exists $\ell_0\in\mN$ such the natural injection $A(K^{p^{-\ell_0}})\hookrightarrow A(K^\perf)$ is surjective (and hence a bijection). 
\label{mprop1}
\end{theor}

For the notion of ampleness, see \cite[par. 2]{Hartshorne-Ample}. A smooth commutative $S$-group scheme $\CA$ as above is called semiabelian if 
each fibre of  $\CA$ is an extension of an abelian variety by a torus (see \cite[I, def. 2.3]{Faltings-Chai-Degeneration} for more details). 

We recall the following fact, which is proven in \cite{Barton-Tensor}: a vector bundle $V$ on $S$ is ample if and only if $\bar\mu_\min(V)>0$. 

\begin{theor}
Suppose that $A$ is an ordinary abelian variety. Then
\begin{itemize}
\item[\rm (a)] $\bar\mu_\min(\omega)\geqslant 0$;
\item[\rm (b)] if there is a closed point $s\in S$ such that $\CA_s$ is an abelian variety of $p$-rank $0$, then $\bar\mu_\min(\omega)>0$. 
\end{itemize}
\label{mprop2}
\end{theor}

\begin{cor} Suppose that $A$ is ordinary and that there is a closed point $s\in S$ such that $\CA_s$ is an abelian variety of $p$-rank $0$.  Then  
\begin{itemize}
\item[\rm (a)]  there exists $\ell_0\in\mN$ such the natural injection $A(K^{p^{-\ell_0}})\hookrightarrow A(K^\perf)$ is surjective;
\item[\rm (b)] the group $A(K^\perf)/\Tr_{K|k}(A)(k)$ is finitely generated.
\end{itemize}
\label{coramp}
\end{cor}  

The notation $\Tr_{K|k}(A)$ refers to the $K|k$-trace of $A$ over $k$. This is an abelian variety 
over $k$, which comes with a morphism $\Tr_{K|k}(A)_K\to A$. See \cite{Conrad-Chow} for the definition. 

Here is an application of Corollary \ref{coramp}. Suppose until the end of the sentence that $A$ is an elliptic curve over $K$ and that $j(A)\not\in k$ (here $j(\cdot)$ is the modular $j$-invariant); then $\Tr_{K|k}(A)=0$, $A$ is ordinary and there is a closed point $s\in S$ such that $\CA_s$ is an elliptic curve of $p$-rank $0$ (i.e. a supersingular elliptic curve); thus $A(K^\perf)$ is finitely generated. This was also proven by D. Ghioca (see \cite{Ghioca-Elliptic}) using a different method.

We list two further applications of Theorems \ref{mprop1} and \ref{mprop2}. 

Let $Y$ be an integral closed subscheme of $B:=A_{\bar K}$. 

Let $C:=\Stab(Y)^\red$, where $\Stab(Y)=\Stab_{B}(Y)$ is the translation stabilizer of $Y$. 
This is the closed subgroup scheme  of $B$, which is characterized uniquely 
by the fact that for any scheme $T$ and any morphism $b: T\to B$, translation by $b$ on the
product $B\times_{\bar K} T$ maps the subscheme $Y\times T$ to itself if and
only if $b$ factors through $\Stab_{B}(Y)$.  Its existence is proven in 
\cite[exp.~VIII, Ex.~6.5~(e)]{SGA3-2}.

The following proposition is a special case of the (unproven) "full" Mordell-Lang conjecture, first 
formulated by Abramovich and Voloch. See \cite{Ghioca-Division} and \cite[Conj. 4.2]{Scanlon-A-positive} for 
a formulation of the conjecture and further references. 

\begin{prop} 
Suppose that $A$ is an ordinary abelian variety.  
Suppose that there is a closed point $s\in S$ such that $\CA_s$ is an abelian variety of $p$-rank $0$.  
Suppose that $\Tr_{\bar K|k}(A)=0$. 
If $Y\cap A(K^\perf)$ is Zariski dense in $Y$ then $Y$ is the translate of an abelian  subvariety of 
$B$ by an point in $B(\bar K)$. 
\label{MLtheor}
\end{prop}

{\bf Proof} (of Proposition \ref{MLtheor}). This is a direct consequence of Corollary \ref{coramp} and of the Mordell-Lang conjecture 
over function fields of positive characteristic; see \cite{Hrushovski-Mordell-Lang}  for the latter.\endProof

Our second application is to a conjecture of a A. Langer and H. Esnault. 
See \cite[Remark 6.3]{Esnault-Langer-On-a-positive} for the latter. The following proposition is a special case of their conjecture. 

\begin{prop} Suppose that $k=\bar\mF_p$.  Suppose that $A$ is an ordinary abelian variety 
and that there is a closed point $s\in S$ such that $\CA_s$ is an abelian variety of $p$-rank $0$. 
Suppose that for all $\ell\geqslant 0$ we are given a point 
$P_\ell\in A^{(p^\ell)}(K)$ and suppose that for all $\ell\geqslant 1$, we have $\Ver^{(p^\ell)}_{A/K}(P_{\ell})=P_{\ell-1}$. 
Then $P_0$ is a torsion point. 
\label{LEconj}
\end{prop}

Here $\Ver^{(p^\ell)}_{A/K}:A^{(p^\ell)}\to A^{(p^{\ell-1})}$ is the Verschiebung morphism. See \cite[$\rm VII_A$, 4.3]{SGA3-1} for the definition. 

{\bf Proof} (of Proposition \ref{LEconj}). By assumption the point $P_0$ is $p^\infty$-divisible in 
$A(K^\perf)$, because $[p]_{A/K}=\Ver_{A/K}\circ\Frob_{A/K}$ and  $\Ver_{A/K}$ is \'etale, because $A$ is ordinary. 
Here $\Frob_{A/K}$ is the relative Frobenius morphism and $[p]_{A/K}$ is the multiplication by $p$ morphism on $A$. Thus the image of 
$P_0$ in $A(K^\perf)/\Tr_{K|k}(A)(k)$ is a torsion point because the group $A(K^\perf)/\Tr_{K|k}(A)(k)$  is finitely generated 
by Corollary \ref{coramp}. Hence $P_0$ is a torsion point because $\Tr_{K|k}(A)(k)$ consists of torsion points.\endProof

{\bf Acknowledgments.} I would like to thank H. Esnault and A. Langer for many interesting conversations on the subject matter of this article.  

\section{Proof of \ref{mprop1}, \ref{mprop2} \& \ref{coramp}}

\subsection{Proof of Theorem \ref{mprop1}}

The idea behind the proof of Theorem \ref{mprop1} comes from an article of M. Kim (see \cite{Kim-Purely-inseparable}). 

\label{pp1}

In this subsection, the assumptions of Theorem \ref{mprop1} hold. So we suppose that $\CA/S$ is semiabelian, that $A$ is a principally polarized abelian variety and that $\omega$ is ample.

If $Z\to W$ is a $W$-scheme and $W$ is a scheme of characteristic $p$, then for any $n\geqslant 0$ we shall write 
$Z^{[n]}\to W$ for the $W$-scheme given by the composition of arrows $$Z\to W\stackrel{F_W^{n}}{\to}W.$$

Now fix $n\geqslant 1$ and suppose that $A(K^{p^{-n}})\backslash A(K^{p^{-n+1}})\not=\emptyset$. 

Fix $P\in A^{(p^n)}(K)\backslash A^{(p^{n-1})}(K)=A(K^{p^{-n}})\backslash A(K^{p^{-n+1}}).$
The point $P$ corresponds to a  commutative diagram of $k$-schemes
\begin{diagram}
& & A\\
& \ruTo^{P} & \dTo\\
\Spec\ K^{[n]} & \rTo_{F_K^n} &\Spec\ K
\end{diagram}
such that the residue field extension $K|\kappa(P(\Spec\ K^{[n]}))$ is of degree $1$ (in other words $P$ is birational onto its image).  In particular, the  
map of $K$-vector spaces $P^*\Omega^1_{A/k}\to\Omega^1_{K^{[n]}/k}$ arising from the diagram is non zero. 

Now recall that there is a canonical exact sequence
$$
0\to \pi^*_K\Omega^1_{K/k}\to \Omega^1_{A/k}\to\Omega^1_{A/K}\to 0.
$$
Furthermore the map $F_K^{n,*}\Omega^1_{K/k}\stackrel{F_K^{n,*}}{\to}\Omega^1_{K^{[n]}/k}$ vanishes. Also, 
we have a canonical identification $\Omega^1_{A/K}=\pi_K^*\omega_K$ (see \cite[chap. 4., Prop. 2]{Bosch-Raynaud-Neron}). Thus the natural surjection 
\mbox{$P^*\Omega^1_{A/k}\to\Omega^1_{K^{[n]}/k}$} gives rise to a non-zero map 
$$\phi_n:F_K^{n,*}\omega_K\to \Omega^1_{K^{[n]}/k}.$$

The next lemma examines the poles of the morphism $\phi_n$.  

We let $E$ be the closed subset, which is the union of the points $s\in S$, such that the fibre $\CA_s$ is not complete. 

\begin{lemma}
The morphism $\phi_n$ extends to 
a morphism of vector bundles $$F_S^{n,*}\omega\to \Omega^1_{S^{[n]}/k}(E).$$
\label{flem}
\end{lemma}
\beginProof (of \ref{flem}).   
First notice 
that there is a natural identification \mbox{$\Omega^1_{S^{[n]}/k}(\log E)=\Omega^1_{{S^{[n]}}/k}(E)$,} because there is a 
sequence of coherent sheaves
$$
0\to\Omega_{{S^{[n]}}/k}\to \Omega^1_{{S^{[n]}}/k}(\log E)\to \CO_E\to 0
$$
where the morphism onto $\CO_E$ is the residue morphism. Here the sheaf $\Omega^1_{{S^{[n]}}/k}(\log E)$ is the 
sheaf of differentials on $S^{[n]}\backslash E$ with logarithmic singularities along $E$. See \cite[Intro.]{Illusie-Reduction} for this result and more details on these notions. 

Now notice that in our proof of Theorem \ref{mprop1}, we may replace $K$ by a finite extension field $K'$ without restriction of generality. 
We may thus suppose that $A$ is endowed with an $m$-level structure for some $m\geqslant 3$. 

We now quote part of one the main results of the book \cite{Faltings-Chai-Degeneration}:
\begin{itemize}
\item[(1)] there exists a regular moduli space $A_{g,m}$ for principally polarized abelian varieties over $k$ endowed with 
an $m$-level structure;
\item[(2)] there exists an open immersion $A_{g,m}\hookrightarrow A^*_{g,m}$, such that 
the (reduced) complement $D:=A^*_{g,m}\backslash A_{g,m}$ is a divisor with normal crossings and 
$A^*_{g,m}$ is regular and proper over $k$;
\item[(3)] the scheme $A^*_{g,m}$ carries a semiabelian scheme $G$ extending the universal abelian scheme $f:Y\to A_{g,m}$;
\item[(4)] there exists a regular and proper $A^*_{g,m}$-scheme $\bar f:\bar Y\to A^*_{g,m}$, which extends $Y$ and such that 
$F:=\bar Y\backslash Y$ is a divisor with normal crossings (over $k$); furthermore 
\item[(5)] on $\bar Y$ there is an exact sequence of locally free sheaves
$$
0\to \bar f^*\Omega^1_{A_{g,m}/k}(\log D)\to\Omega^1_{Y/k}(\log F)\to\Omega^1_{Y/A_{g,m}}(\log F/D)\to 0,
$$
which extends the usual sequence of locally free sheaves
$$
0\to f^*\Omega^1_{A_{g,m}/k}\to\Omega^1_{Y/k}\to\Omega^1_{Y/A_{g,m}}\to 0
$$
on $A_{g,m}$. Furthermore there is an isomorphism  $\Omega^1_{Y/A_{g,m}}(\log F/D)\simeq\bar f^*\omega_G$. Here $\omega_G:={\rm Lie}(G)^\vee$ is the tangent bundle (relative to $A_{g,m}^*$) of $G$ restricted to $A_{g,m}^*$ via the unit section. 
\end{itemize}
See \cite[chap. VI, th. 1.1]{Faltings-Chai-Degeneration} for the proof. 

The datum of $A/K$ and its level structure induces a morphism $\phi:K\to A_{g,m}$, such that 
$\phi^*Y\simeq A$, where the isomorphism respects the level structures. Call $\lambda:A\to Y$ 
the corresponding morphism over $k$. Let $\bar\phi:S\to A^*_{g,m}$ be the morphism obtained 
from $\phi$ via the valuative criterion of properness. By the unicity of semiabelian models 
(see \cite[chap. I, th. I.9]{Faltings-Chai-Degeneration}), we have a natural isomorphism $\bar\phi^*G\simeq\CA$ and thus we have a set-theoretic equality $\bar\phi^{-1}(D)=E$ 
and an isomorphism $\bar\phi^*\omega_G=\omega$. 
Let also $\bar P$ be the morphism $S^{[n]}\to\bar Y$ obtained 
from $\lambda\circ P$ via the valuative criterion of properness. By construction we now get an arrow 
$$
\bar P^*\Omega^1_{\bar Y/k}(\log F)\to\Omega^1_{S^{[n]}/k}(\log E)
$$
and since 
the induced arrow 
$$
\bar P^*\bar f^*\Omega^1_{A_{g,m}/k}(\log D)=F^{n,*}_S\circ\bar\phi^*(\Omega^1_{A_{g,m}/k}(\log D))\to 
\Omega^1_{S^{[n]}/k}(\log E)
$$
vanishes (because it vanishes generically), we get an arrow 
$$
\bar P^*\Omega^1_{\bar Y/A^*_{g,m}}(\log F/D)=F^{n,*}_S\circ\bar\phi^*\omega_G=
F^{n,*}_S\omega\to \Omega^1_{S^{[n]}/k}(\log E)=\Omega^1_{S^{[n]}/k}(E),
$$
which is what we sought.\endProof

To conclude the proof of Proposition \ref{mprop1}, choose $l_0$ large enough so that 
$$\mu_\min(F^{l,*}_S(\omega))>\mu(\Omega^1_{S/k}(E))$$
for all $l>l_0$. Such an $l_0$ exists because $\bar\mu_\min(\omega)>0$. 
Now notice that since $k$ is a perfect field, we have $\Omega^1_{S/k}(E)\simeq \Omega^1_{S^{[n]}/k}(E)$. 
We see that we thus have 
$$\Hom(F^{l,*}_S(\omega),\Omega^1_{S^{[n]}/k}(E))=0$$ for all $l>l_0$ and thus by Lemma \ref{flem} 
we must have $n<l_0+1$. Thus  we have 
$$
A(K^{(p^{-l})})=A(K^{(p^{-l+1})})
$$
for all $l\geqslant l_0$. 

{\bf Remark}. The fact that $\Hom(F^{l,*}_S(\omega),\Omega^1_{S^{[n]}/k}(E))\simeq \Hom(F^{l,*}_S(\omega),\Omega^1_{S/k}(E))$ vanishes for large $l$ can also be proven without appealing 
to the Harder-Narasimhan filtration. Indeed the vector bundle $\omega$ is also cohomologically $p$-ample (see \cite[Rem. 6), p. 91]{Migliorini-Some}) and thus there is an $l_0\in\mN$ such that for all 
$l>l_0$  
\begin{eqnarray*}
\Hom(F^{l,*}_S(\omega),\Omega^1_{S/k}(E))&=&H^0(S,F^{l,*}_S(\omega)^\vee\otimes \Omega^1_{S/k}(E))\\
&\stackrel{{\rm Serre\ duality}}{=}&
H^1(S,F^{l,*}_S(\omega)\otimes \Omega^1_{S/k}(E)^\vee\otimes\Omega^{1}_{S/k})^\vee\\
&=&H^1(S,F^{l,*}_S(\omega)\otimes\CO(-E))^\vee=0.
\end{eqnarray*}

\subsection{Proof of Theorem \ref{mprop2}}

\label{pp2}
In this subsection, we suppose that the assumptions of Theorem \ref{mprop2} hold. So we suppose that 
$A$ is an ordinary abelian variety. 

Notice first that  for any $n\geqslant 0$, the Hodge bundle of $\CA^{(p^n)}$ is $F_S^{n,*}\omega$. 
Hence, in proving Proposition \ref{mprop2}, we may assume without restriction of generality that $\omega$ has a strongly semistable Harder-Narasimhan filtration. 

Let $V:=\omega/\omega_{\hn(\omega)-1}$. Notice that for any $n\geqslant 0$, we have a (composition of) Verschiebung(s) map(s)
$
\omega{\to} F_S^{n,*}\omega.
$
Composing this with the natural quotient map, we get a map 
\begin{equation}
\phi:\omega\stackrel{\Ver^{(p^n),*}_{\CA}}{\to} F_S^{n,*}V
\label{sfeq}
\end{equation}
The map $\phi$ is generically surjective, because by the assumption of ordinariness the map \mbox{$\omega\stackrel{\Ver^{(p^n),*}}{\to} F_S^{n,*}\omega$} is generically an isomorphism. 

We now prove (a). The proof is by contradiction. Suppose that $\bar\mu_\min(\omega):=\mu(V)<0$. This implies that when $n\to\infty$, we have 
$\mu(F_S^{n,*}V)\to-\infty$. Hence if $n$ is sufficiently large, we have $\Hom(\omega,F_S^{n,*}V)=0$, which contradicts 
the surjectivity of the map in \refeq{sfeq}. 

We turn to the proof of (b). Again the proof is by contradiction. So suppose that $\bar\mu_\min(\omega)\leqslant 0$. 
By (a), we know that we then actually have \mbox{$\bar\mu_\min(\omega)=0=\mu(V)$} and $V\not=0$. If $\bar\mu_\max(\omega)>0$ then the map $\omega_{1}\to F_S^{n,*}V$ obtained by composing $\phi$ with the inclusion $\omega_1\hookrightarrow\omega$ must vanish, because 
$$\mu(\omega_1)>\mu(F_S^{n,*}V)=p^n\cdot\mu(V)=0.$$ Hence we obtain a map $\omega/\omega_1\to F_S^{n,*}V$. 
Repeating this reasoning for $\omega/\omega_1$ and applying induction  we finally get 
a map 
$$
\lambda:V\to F_S^{n,*}V.
$$
The map $\lambda$ is generically surjective and thus globally injective, since its target and source are locally free sheaves of the same generic rank. Let $T$ be the cokernel of $\lambda$ (which is a torsion sheaf). We then have 
$$
\deg(V)+\deg(T)=0+\deg(T)=\deg(F_S^{n,*}V)=0
$$
and thus $T=0$. This shows that $\lambda$ is a (global) isomorphism. 
In particular, the map $\phi$ is surjective. Thus the map
$$
\phi_s:\omega_s\stackrel{\Ver^{(p^n),*}_{\CA_s}}{\to} F_s^{n,*}V_s
$$
is surjective and thus non-vanishing. This contradicts the hypothesis on the $p$-rank at $s$.

\subsection{Proof of Corollary \ref{coramp}}

In this subsection, we suppose that the assumptions of Corollary \ref{coramp} are satisfied. So we suppose that $A$ is ordinary and that there is a closed point $s\in S$ such that $\CA_s$ is an abelian variety of $p$-rank $0$.

We first prove (a). First we may suppose without restriction of generality that $A$ is principally polarized. 
This follows from the fact that the abelian variety $(A\times_K A^\vee)^4$ carries a principal polarization ("Zarhin's trick" - see 
\cite[Rem. 16.12, p. 136]{Milne-Abelian}) and from the fact that the abelian variety $(A\times_K A^\vee)^4$ also 
satisfies the assumptions of Corollary \ref{coramp}. Furthermore, we may without restriction of generality replace $S$ by a finite extension $S'$. 
Thus, by Grothendieck's semiabelian reduction theorem (see \cite[IX]{SGA7.1}) we may assume that $\CA$ is semiabelian. Statement (a) 
then follows from Theorems \ref{mprop1} and \ref{mprop2}. 

We now turn to statement (b). 
Let $\tau_{K|k}:\Tr_{K|k}(A)_K\to A$ be the $K|k$-trace morphism. Notice that for any $\ell\geqslant 0$, we have a natural identification of 
$k$-group schemes $\Tr_{K|k}(A)^{(p^\ell)}\simeq\Tr_{K|k}(A^{(p^\ell)})$, because the extension $K/K^p$ is primary 
and $k$ is perfect (see \cite[Th. 6.4 (3)]{Conrad-Chow}). 
Thus, if $\ell_0\in\mN$ is the number appearing in (a), we have identifications
\begin{eqnarray*}
&&A(K^\perf)/\Tr_{K|k}(A)(k)=A(K^{-\ell_0})/\Tr_{K|k}(A)(k)=A(K^{-\ell_0})/\Tr_{K|k}(A)(k^{-\ell_0})\\
&&=A^{(p^{\ell_0})}(K)/\Tr_{K|k}(A)^{(p^{\ell_0})}(k)=A^{(p^{\ell_0})}(K)/\Tr_{K|k}(A^{(p^{\ell_0})})(k)
\end{eqnarray*}
and the group appearing after the last equality is finitely generated by the Lang-N\'eron theorem.

\end{document}